\documentclass[12pt]{article}
\usepackage{amssymb}
\usepackage{amsthm}
\usepackage{lineno}
\usepackage{amsfonts}
\usepackage{amsmath}

\newtheorem{thm}{Theorem}

\newtheorem{proposition}{Proposition}

\begin{document}
\title{Zeckendorf arithmetic for Lucas numbers}
\author{\textbf{Rachid CHERGUI} \thanks{%
The research is partially supported by ATN laboratory of USTHB University.%
} \\
{\small USTHB, Faculty of mathematics }\\
[-0.8ex] {\small P.B. 32, El Alia, 16111, Bab Ezzouar, Algeria}\\
[-0.8ex] {\small \texttt{rchergui@usthb.dz 
}}
}
\date{ {\small Mathematics Subject Classifications: 05A10,
05A17, 11B39, 11B50.}}
\maketitle
\begin{abstract}
In this article we will be dedicated some algorithms of addition, subtraction, multiplication and division of two positive integers using Zeckendorf form. Such results find application in coding theory.
\end{abstract}
\maketitle
\section{Introduction}
There are few previous work in this area. Graham, Knuth and Patashnik \cite{6} discussing the addition of 1 in the Zeckendorf representation, but have
not talked about the actual arithmetic. Fliponi \cite{4} did for addition and
multiplication, and Freitag philips \cite{2} for the subtraction and division
\cite{3}. Thus, no previous work has discussed arithmetic as a coherent
whole, covering all major operations, including multiplication and division.
All these algorithms have been implemented and tested on a computer. Most
algorithms are developed by analogy with conventional arithmetic methods. For example, multiplication is carried out by adding appropriate multiples
of the multiplicand, depending on the selected bit pattern of the multiplier. The division will use a sequence of test subtraction, as in the normal long
division.
\section{Zeckendorf theorems for Lucas numbers } 
Lucas numbers are defined by the recursion formula:
\[
 \left \{
\begin{array}{ll}
L_{n} =L_{n-1}+L_{n-2} ,n\geq2 \\
L_{0}=2,L_{1}=1.
\end{array}
 \right .\]
and for all $n\geq0$, we have the well-known $L_{n} =F_{n+1}+F_{n-1}$ where $F_{n}$ is the nth Fibonacci number.

\subsection{ Zeckendorf decompositon method}
To decompose an integer $x$ of the form Zeckendorf $x=\sum_{r=0}^{\propto}\alpha_{r}L_{r} $  just follow the following steps:\\
1.Find the greater Lucas number $ L_{r}\leqslant x$.  \\
2.Do subtraction $X=x-L_{r}$, assign a 1 to $e_{r} $ and keep this coefficient.\\
3.Assign $X$ to $x$ and repeat steps 1 et 2 until $X$ have a zero.\\
4.Assign  of  0 to  $e_{i}$ where $0 \leqslant i \leqslant r$ and $e_{i}\neq1 $.\\
The result of this decomposition is a vector of $r$  elements that contains
the coefficients $e_{i}$ decomposition.
Example decomposition of $ 50$, this table shows the performance:
\begin{center}
\begin{tabular}{|l|l|l|l|l|l|l|l|l|l|l|}
\hline
Lucas sequence& 2 & 1 & 3 & 4 & 7 & 11 & 18 & 29 & 47 & 76 \\
\hline
Vector $e_{r}$ & 0 & 0 & 1 & 0 & 0 & 0 & 0 & 0 & 1 & 0 \\
\hline
\end{tabular}
\end{center}
 So $50_{L} = 001000001$
 \begin{thm}
Let $n$ be a positive integer satisfying $0\leqslant n\leqslant L_{k} $
for $ k\geqslant 1$. Then \\
$n=\sum_{i=0}^{k-1}\alpha_{i}L_{i} $ where  $\alpha_{i} \in\lbrace0,1\rbrace$ such that \\
\[
 \left \{
\begin{array}{ll}
i)\alpha_{i}\alpha_{i+1}=0, for (i\geqslant0) \\
ii)\alpha_{0}\alpha_{2}=0.
\end{array}
 \right .
\]
this representation is unique.
\end{thm}
Proof: \cite{1}

\section{Addition}
\noindent
We take two positive integers a and b written in the form of Zeckendorf,
obtainable form of a + b Zeckendorf repeating adding, at the same time,
numbers of Lucas occupant in one of two numbers , say b , to another
number a. This gives an initial amount for which figures are $d_{i}\in \left\{ 0,1\right\} $ , where each $d_{i}$ is $L_{i}$ its number of Lucas. For $d_{i}= 2$ does not exist because $n=1 \rightarrow 2L_{1}=L_{0}$, we replace $020$ by $001$ and $n\geq2\rightarrow 2L_{n} = L_{n+1} + L_{n-2}$, we replace $00200$ by $01001$. In way is equivalent model x 2 y z figures transforms to\\
(1 + x )0y (1 + z). This rule does not apply to terms with a weight of 1, which is covered by the special case below. If the combination $011$ exists in the vector $e_{r}$, we will substitute it by $001$. This step must be performed by scanning left to right through the performance. Here is a table that summarizes all possible cases of the addition in the representation of Zeckendorf:
\newpage
\begin{table}[h]
\begin{center}
\begin{tabular}{|ll@{\extracolsep{.50cm}}crcccccc|}
\hline
 & &  &  & & \multicolumn{4}{c}{}& \\
& \textbf{Addition} & &  Lucas weight  & & $L_{i+1}$ & $L_{i}$ & $L_{i-1}$ & $L_{i-2}$ & \\\cline{5-8}
& Consecutive 1 & & & & $x$ & $y$ & $1$ & $1$ & \\
&  & & becomes & & $x$ & $y + 1$ & $0$ & $0$ & \\
& & & & & & & & & \\
&Eliminate a 2& & here $x\geq2$ & & $w$ & $x$ & $y$ & $z$ & \\
 && & becomes & & $w + 1$ & $x -2$ & $y$ & $z + 1$ & \\\hline
& &  &  & & \multicolumn{3}{c}{}& & \\
& \textbf{Add,right bits } & & & &  $L_{2}$ & $L_{1}$ & $L_{0}$ & &\\\cline{5-7}
&$d_{2}\geq 2$ & & here $x\geq2$ & & & $x$ & $0$ & & \\
& & & becomes & & & $0$ & 1 & & \\
& & & & & & & & & \\
&$d_{2}\geq2$ & & & &  & $0$ & $x$ & & \\
& & & becomes & & $ 1$ & $1$ & 0 & & \\
 \\\hline
\end{tabular}
\end{center}
\caption{Adjustments and corrections in  addition}
\end{table}
This table shows the two additions examples 33+19 and 12+19 in zeckendorf representation:
\begin{table}[h]
\begin{center}
\begin{tabular}{|ll@{\extracolsep{.50cm}}crccccccccc|}
\hline
&a & & &1& 0& 0& 0& 1& 0& 0& 0& = 33 \\
&b & & & &1& 0& 0& 0& 0& 1& 0&  = 19 \\
&initial sum & & & 1& 1& 0& 0& 1& 0& 1& 0&  = 52 \\
& consecutive 1 & &1& 0& 0& 0& 0& 1& 0& 1& 0&  = 52 \\\cline{4-12}
&becomes & & 1& 0& 0& 0& 0& 1& 0& 1& 0&  = 52 \\
&\multicolumn{11}{l}{check $ 33 +19 = 52$} &\\\hline
\end{tabular}
\end{center}
\caption{Example  of addition $(33+19)$}
\end{table}
\begin{table}[h]
\begin{center}
\begin{tabular}{|ll@{\extracolsep{.50cm}}crccccccccc|}
\hline
&a& &&&1& 0& 0& 0& 1& 0& &= 12\\
&b&&& 1& 0& 0& 0& 0& 1& 0& &= 19\\
&initial sum &&& 1& 1& 0& 0& 0& 2& 0& &= 31\\
&     &&& 1& 1& 0& 0& 0& 0& 1&& = 31\\
&consecutives 1& && 1& 1& 0& 0& 0& 0& 1& & = 31\\
&     & & 1& 0& 0& 0& 0& 0& 0& 1&& = 31\\\cline{4-11}
&becomes&& 1& 0& 0& 0& 0& 0& 0& 1&& = 31\\
&\multicolumn{11}{l}{check $ 12+ 19 = 31$}&\\\hline
\end{tabular}
\end{center}
\caption{Example of addition $(12+19)$}
\end{table}
\vfil\eject
\section{Soustraction}
For subtraction, $ X-Y = Z$, where $ X > Y $ and Z the difference. We start by
subtracting all figures $x_{i}-y_{i} =z_{i}$, where $z_{i}\in \left\{ 0,1.-1\right\} $ Values 0 and 1 have no problem, as they are valid representation in Zeckendorf.
Where $z_{i} = -1 $ is the most difficult. If in this case, go to the next bit
1 and is written in the Fibonacci rule  $100\rightarrow  011$ and write bit 1 rightmost
pairs 1 is repeated until the bit 1 bit coincides with the -1 in the same position and eliminates replacing the box by 0 then 1 consecutive passes.
\begin{table}[h]
\begin{center}
\begin{tabular}{|l@{\extracolsep{.50cm}}crccccccc|}
\hline
  &  &  & & \multicolumn{5}{c}{}& \\
\textbf{Soustraction} & & Lucas weights  & & $L_{i+2}$ & $L_{i+1}$ & $L_{i}$ & $L_{i-1}$ & $L_{i-2}$ & \\\cline{5-9}
eliminate -1 & & & & 1 & 0 & 0 & 0 & -1 & \\
 & &  & & 0 & 1 & 1 & 0 & -1 & \\
 & & becomes & & 0 & 1 & 0 & 1 & 0 & \\\hline
\end{tabular}
\end{center}
\caption{ Adjustments and corrections in  subtraction}
\end{table}\\
This table shows the example 42-32 in Zeckendorf representation.As elsewhere,different rewriting rules can be written in any order.
\begin{table}[h]
\begin{center}
\begin{tabular}{|ll@{\extracolsep{.50cm}}crcccccccc|}
\hline
&a& &1 &0 &1 &0 &0 &0 &0 &1& = 42\\
&b& &1 &0 &0 &0 &0 &1 &0 &0 &= 32\\\hline
&subtract bit by bit&&&& 1& 0& 0&  -1& 0& 1& = 10\\
&rewrite 1000& &&&&1& 1& -1& 0& 1& = 10\\
&rewrite 0110, cancelling -1 &&&&& 1& 0& 0& 1& 1& = 10\\
&consecutive 1&&&&& 1& 0& 1& 0& 0& = 10\\
&becomes &&&&& 1& 0& 1& 0& 0& = 10\\\hline
\end{tabular}
\end{center}
\caption{Example of subtraction ($42-32)$}
\end{table}
\newpage
\section{Multiplication}
Using the following results (propositions 1,2,3,4) and section 3 above, one can derive a multiplication method of integers in Zackendorf  representaion. \\
\begin{proposition}
If $n\geq3$, then
\[
 L_{k}L_{k+n}=\left \{
\begin{array}{ll}
 F_{n-1}+F_{n+1}+F_{2k+n\pm1} & (k \hspace{.25cm} even), \\
 F_{n-2}+F_{n+1}+F_{2k+n+1}+\sum_{j=1}^{k-2}F_{2j+n+2} & (k\geq3, odd).
\end{array}
 \right .
\]
\end{proposition}
\begin{proposition}
If $n\geq5$, then
\[
 2L_{k}L_{k+n}=\left \{
\begin{array}{ll}
 F_{n\pm3}+F_{2k+n\pm3} & (k\geq4, even), \\
 F_{n-4}+F_{2k+n+3}+\sum_{j=1}^{3}F_{2j+n-3}+\sum_{j=1}^{k-4}F_{2j+n+4} & (k\geq5, odd).
\end{array}
 \right .
\]
\end{proposition}
\begin{proposition}
If $n\geq5$, then
\[
 3L_{k}L_{k+n}=\left \{
\begin{array}{ll}
 \sum_{j=1}^{4}\left(F_{2j+n-5}+F_{2j+2k+n-5}\right) & (k\geq4, even), \\
 F_{n-4}+F_{n+3}+\sum_{j=1}^{3}F_{2j+2k+n-3}+\sum_{j=1}^{k-4}F_{2j+n+4} & (k\geq5, odd).
\end{array}
 \right .
\]
\end{proposition}
\begin{proposition}
If $n\geq6$, then
\[
 4L_{k}L_{k+n}=\left \{
\begin{array}{ll}
 \sum_{j=1}^{4}\left(F_{3j+n-8}+F_{3j+2k+n-8}\right) & (k\geq6, even), \\
 F_{n-4}+F_{n-2}+F_{n+1}+\sum_{j=1}^{3}F_{3j+2k+n-5}+\sum_{j=1}^{k-5}F_{2j+n+4} & (k\geq5, odd).
\end{array}
 \right .
\]
\end{proposition}
Proofs: \cite{5}\\

This example shows how to compute  $17 \times 10$ in Zeckendorf representation:
\begin{table}[h]
\begin{center}
\begin{tabular}{|l@{\extracolsep{.50cm}}cccccccccccc|}
\hline
a&&&&&  &1 &0 &1 &0 &0 &1 &=17 \\
b&&&&& & &1& 0& 1& 0& 0&=10\\\hline
\multicolumn{13}{|l|}{\textbf{Multiple of Luca 17}}\\
multiple $L_{1}$ & &&&  & & 1&0 &1 &0 &0 &1 &=17\\
multiple $L_{2}$ & &&1& 0 &0 & 0&0 &1 &0 &0 &0 &=51\\
multiple $L_{3}$ & && 1 &0 &1 &0 &0 &0 &1 &0 &0&=68\\
multiple $L_{4}$ & & 1& 0& 1& 0& 1& 0& 0& 1& 0& 0&=119\\
multiple $L_{5}$ &1&0& 1&  0& 0&1 &0 &1 &0 &0 &1 &=187\\\hline
\multicolumn{13}{|l|}{\textbf{Accumulate appropriate multiples }}\\
Add multiple of $L_{2}$& &&1& 0 &0 & 0&0 &1 &0 &0 &0 &=51\\

Add multiple of $L_{4}$ & & 1& 0& 1& 0& 1& 0& 0& 1& 0& 0&=119\\
$L_{2}+L_{4}$=& & 1& 1& 1& 0& 1& 0& 1& 1& 0& 0& =170\\
Eliminate 1 consecutive  & 1& 0& 0& 1& 0& 1& 1& 0& 0& 0& 0& =170\\
Eliminate 1 consecutive & 1& 0& 0& 1& 1& 0& 0& 0& 0& 0& 0& =170\\
becomes=& 1& 0& 1& 0& 0& 0& 0& 0& 0& 0& 0& =170\\\hline

\end{tabular}
\end{center}
\caption{Example of Zeckendorf multiplication $(17\times10)$}
\end{table}
\vfil\eject
\section{Division}
Using the following proposition 5 and section 4 above, one can derive a division method of integers in Zackendorf  representation.\\
\begin{proposition}

First,for $k=4m$ and $n$ odd, we obtain
$$\frac{F_{kn}}{F_n}=\sum_{r=1}^{m}\left(L_{(k-4r+3)n}+L_{(k-4r+1)n}\right) ,$$
and thus
$$\frac{F_{kn}}{F_n}=S_{k,n},$$
say, where
$$S_{k,n}=\sum_{r=0}^{[k/4]-1}\left(F_{(k-4r-1)n+1}+\left(\sum_{s=1}^{n-2}F_{(k-4r-1)n-2s}
\right)+F_{(k-4r-3)n+1}+F_{(k-4r-3)n-2}\right)   .$$
We similarly work through the other cases, where $n$ is odd and $k\equiv1,2$ and $3\mod4$. In each case, the "most significant" part of the Zeckendorf form is $S_{k,n}$. The precise Zackendorf form is
$$\frac{F_{kn}}{F_n}=S_{k,n}+e_{k,n},$$
where the least significant part of the Zeckendorf sum is
\[
 e_{k,n}=\left \{
\begin{array}{ll}
0, &  k\equiv0\mod4), \\
F_{2}, & k\equiv1\mod4), \\
F_{n+1}+F_{n-1}, & k\equiv2\mod4), \\
F_{2n+1}+\sum_{r=1}^{n-1}F_{2n-2r}, & k\equiv3\mod4). \\
\end{array}
 \right .
\]
\end{proposition}
Proof:\cite{3} \\

This example shows how to compute $250\div17$ in Zeckendorf representation:

\begin{table}[h]
\begin{center}
\begin{tabular}{|l@{\extracolsep{.50cm}}cccccccccccccc|}
\hline
a && &1 &0 &1 &0 &1 &0 &0 &0 &1 &0 &0 &=250\\
b & &&&&&&&1 &0 &1 &0 &0 &1 &=17\\\hline
\multicolumn{15}{|l|}{\textbf{Make lucas Multiples of divisor }}\\
 multiple $L_{1}$&&&&&&&& 1 &0 &1 &0 &0 &1 &=17\\
 multiple $L_{2}$ &&&&&1& 0 &0 &0 &0 &1 &0 &0 &0 &=51\\
 multiple $L_{3}$ & &&&&1 &0 &1 &0 &0 &0 &1 &0 &0 &=68\\
 multiple $L_{4}$ & &&&1 &0 &1 &0 &1 &0 &0 &1 &0 &0 &=119\\
 multiple $L_{5}$ &&& 1 &0 &1 &0 &0 &1 &0 &1 &0 &0 &1 &=187\\
 multiple $L_{6}$ && 1 &0 &1 &0 &1 &0 &0 &0 &0 &0 &0 &1&=306\\\hline
\multicolumn{15}{|l|}{\textbf{Trial subtraction }}\\
$L5$ residue=& &&&& 1&0 &0 &1&0 &1 &0 &1 &0 &=63\\
$L2$ residue=&&&&&&&& 1&0 &0 &0 &1 &0 &=12\\\cline{10-15}
quotient=&&&&&&& &1 &0 &1 &0 &0 &1 &=17\\
remainder =&&&&&&&& 1&0 &0 &0 &1 &0 &=12\\\hline

\end{tabular}\\
\end{center}
\caption{Example of Zeckendorf division $(250 \div17)$ }
\end{table}
\newpage

\section{Conclusion}
Although we have highlighted the main arithmetic operations on integers
Zeckendorf, this arithmetic should not stay more than a curiosity. In future
research, we plan to study the applications of our results to other areas of
mathematics such as error correcting codes.

\medskip


\end{document}